\documentclass[12pt]{amsart}
\usepackage{geometry} % see geometry.pdf on how to lay out the page. There's lots.
\usepackage{amsfonts}
\usepackage{cite}
\usepackage{graphicx}
\usepackage{amssymb,latexsym,amsmath,mathrsfs} 
\usepackage{listings}

\title{A GPU implementation of the Discontinuous Galerkin method for simulation of diffusion in brain tissue}
\author{Daniel Cervantes \and Miguel angel Moreles \and Joaquin Pe\~na \and Alonso Ramirez-Manzanares}
\address{D. Cervantes \\
               INFOTEC \\
               Av. San Fernando 37, Col. Toriello Guerra, Tlalpan\\
                 Ciudad de Mexico 14050, Mexico\\
               M. A. Moreles, J. Pe\~{n}a, A. Ramirez-Manzanares \\
                 \email{moreles@cimat.mx}\\
              Centro de Investigaci\'{o}n en Matem\'{a}ticas\\ Jalisco s/n, Valenciana\\
Guanajuato, GTO 36240,Mexico           
}

\date{} % delete this line to display the current date

\begin{document}
\begin{abstract}
In this work we develop a methodology to approximate the covariance matrix associated to the simulation of water diffusion inside the brain tissue. The computation is based on an implementation of the Discontinuous Galerkin method of the diffusion equation, in accord with the physical phenomenon. The implementation in in parallel using GPUs in the CUDA language. Numerical results are presented in 2D problems.
\end{abstract}

\maketitle

\tableofcontents

\section{Introduction}

The study of the diffusion phenomenon in porous media is an active line of research with many applications in the applied sciences. The lack of analytical solutions for arbitrary shaped domains requires the use of numerical solvers to compute reliable solutions. The analytical solutions are restricted to simple--shaped geometries as spheres, parallel walls and cylinders \cite{article-Neuman1974}. However, real--life applications require to compute estimations of the diffusion on complex domains.

An important research area in the medical sciences, and the focus of this work, is the simulation of water diffusion inside the brain tissue (white and grey matter). The problem is to characterize the hydrogen molecular displacement due to brownian motion. Literature on the subject is vast, let us present a short review.

The probability of the displacement $x$ (in micrometers) of an ensemble of water molecules in a given time $\Delta$  (in miliseconds) is summarized by the Ensemble Average Propagator (EAP) $P(x,\Delta)$ \cite{price1997}.  The Diffusion Weighted modality of the Magnetic Resonance allows to capture data to infer the molecular displacement on \textit{in vivo} and \textit{ex vivo} experiments \cite{IntroductionCDMRI,article-Alexander2009}. The quantification of molecular displacement of water molecules on \textit{in vivo} patientes allows to infer properties of the microstructure of brain tissue \cite{article-VanGelderen1994,article-CompartmentModelsPanagiotaki2012}. The information above is used to detect tissue disruption associated with deseases: for instance, a premature reduction of the cellular volume is an indicator of neuron death.  The numerical simulation of the water diffusion allows to generate useful data for the validation of model fitting \cite{Fieremans2010MonteCS}, generation of novel theories about the diffusion properties with analytical representations \cite{article-timeDependenceBurcaw2015}, model improvement \cite{Gilani2017AnIM}, etc. In this problem the domain is composed of cell bodies (soma, axons, neurites, glial cells, vessels, etc). Despite the fact that some of the cell structures can be approximated with simple geometries, in general the domains are much more geometrically complex, in particular, the extracellular spaces presents arbitrary shapes (similar to a "gruyere--cheese shape") hence numerical approximations are required to simulate the molecular diffusivity \cite{article-MattHall2009}.

Numerical solvers to compute the EAP of the hydrogen displacement are varied. For methods based on  Monte Carlo diffusion simulators see \cite{article-MattHall2009,Fieremans2010MonteCS,Gilani2017AnIM,article-MattHall2009} and for partial differential--equation based solvers \cite{Nguyen2014AFE,Moroney_MRM24501}. However, the need to produce massive simulation data for: a) simulating experiments for different machine parameters (magnitude of magnetic gradientes, experimental times, etc) \cite{Ferizi2015},  b) validating complex models \cite{article-CompartmentModelsPanagiotaki2012}, c) training automatic learning algorithms \cite{Nedjati_neuroimage_2017}, etc., requires  to produce accurate simulations in optimal computational times.  Nowadays, the diffusion simulators in the state--of--the--art \cite{article-MattHall2009} requiere from minutes to days to estimate the diffusion phenomena on complex domains  \cite{ramirez-manzanaresISMRM2018}.

This leads to the aim of this work, the approximation of the extracellular diffusivity profile on a disordered medium. On a disordered model of cylindrical brain axons, with the diameters computed from a Gamma distribution \cite{article-DistAboitiz1992, article-MattHall2009} the whole process can be characterized by the corresponding 3D covariance matrix $\Sigma$. The $\Sigma$ eigenvalues (matched with the corresponding eigenvectors) indicate the magnitude and orientational dependency, such that, it is possible to infer extracellular microstructure features as: the main orientation of the axon bundles, the percentage of the volume occupied by neurons (intra cellular signal fraction), the amount of diffusion anisotropy of the tissue (fractional anisotropy), among other descriptors. In the DWMR medical literature, this covariance matrix is computed by the diffusion tensor (DT) from the MR signal \cite{ Pierpaoli1996,article-CompartmentModelsPanagiotaki2012}. Those descriptors computed from $\Sigma$ have been correlated with several brain damages and diseases \cite{ DiffusionMRI }.

Consequently, the approximation of the referred covariance matrix is of great interest. We shall introduce a numerical methodology for approximation. As a first step we develop a 2D version of the problem.

The covariance matrix is obtained from a gaussian density, formed by averaging densities which result from the solution of diffusion equations.
Our main contribution is to solve these diffusion equations using and ah hoc implementation of the Discontinuous Galerkin (DG) method, see  \cite{cockburn2012discontinuous} for a thorough discussion. 

Our implementation takes advantage of  the underlying physical and geometric properties of the problem as posed by the clinicians. In practice, the so called substrates are squared pixel domains allowing for a uniform mesh. Also there is an assumption of no diffusivity between the axons and the extracellular region. The common approximation solves the diffusion equation only in the extracellular region imposing a zero Neumann condition. An alternative is to use the numerical fluxes, a main feature of the DG method,  to propose an interaction between the axons and the intercellular region. With this interaction, we are solving the diffusion equation in the whole domain. Thus null computations are carried out in the axons. This, apparently redundant strategy, free us of boundary handling, a computationally expensive task in Galerkin type methods. But more importantly, it allows for the full strength of the DG method, to carry out the time update of the solution in parallel for all elements. A small ODE problem is solved in each element with no communication. Consequently, a GPU implementation is appropriate.

The outline is as follows.

In Section 2, the Initial Boundary Value Problem (IVBP) of the diffusion equation associated to the phenomenon is introduced. The a substrate for study is described for a simulated \textit{ex vivo} experiment.

The basics of the DG method are presented in Section 3. Therein, the physical and geometric properties of the case study are used to tune our DG-GPU implementation. 

In Section 4 a scheme to approximate the diffusion encoding covariance matrix is introduced. Performance is illustrated with free diffusion and the case study
associated to an \textit{ex vivo} experiment. Conclusions and a brief discussion on future work close our exposition.

\section{Water diffusion inside the brain tissue}

In this section we introduce a substrate with such a realistic shape and properties. First, we discuss the Initial Boundary Value Problem (IBVP) for the diffusion equation to be used throughout. 

\subsection{The diffusion problem}

In practice, one models a substrate occupying a domain, which is a medium comprised of two regions, the axons and the extracellular complement. The former is regarded as non diffusive and the latter a region with constant diffusion. It is assumed that the boundary between both regions is reflecting. Initial pulses are prescribed in the extracellular region, far away from the outer boundary which is modelled as a perfect absorber.

Let $\Omega$ be the domain occupied by the substrate. This domain is the union of two intertwined adjacent regions, $\Omega_a$ and $\Omega_e$. These are respectively, the the axon and extracellular regions.

The IVBP consists on finding $u$ that solves the diffusion equation
\begin{equation}
\frac{\partial u}{\partial t}=\nabla\cdot \left(k(x)\nabla u\right),\quad (x,t)\in\Omega_e\times(0,T),
\label{diff_eq}
\end{equation}
given Cauchy data
\begin{equation}
u(x,0)=u_0(x),\quad x\in\Omega_e,
\label{C_data}
\end{equation}
and boundary values
\begin{equation}
\frac{\partial u}{\partial\mathbf{n}}=0,\quad (x,t)\in\partial\Omega_e\times(0,T),
\label{Br_data}
\end{equation}
\begin{equation}
u(x,t)=0,\quad (x,t)\in\partial\Omega\times(0,T).
\label{Bp_data}
\end{equation}

Here $\mathbf{n}$ is the outer unit normal to $\Omega_e$. Notice that $\partial\Omega\subset\partial\Omega_e$.

The Neumann boundary condition (\ref{Br_data}) corresponds to a reflecting boundary, whereas the Dirichlet boundary condition (\ref{Bp_data}) to
that of a perfect barrier.

The discontinuous diffusion is given by
\begin{equation}
k(x)=\left\lbrace
\begin{array}{l}
0,\quad x\in\Omega_a \\
 \\
 k_0,\quad x\in\Omega_e
\end{array}
\right.
\label{k_coeff}
\end{equation} 
for a positive constant $k_0$.

\subsection{A case study}

\begin{figure}[t]    
\centering
\includegraphics[width=300pt]{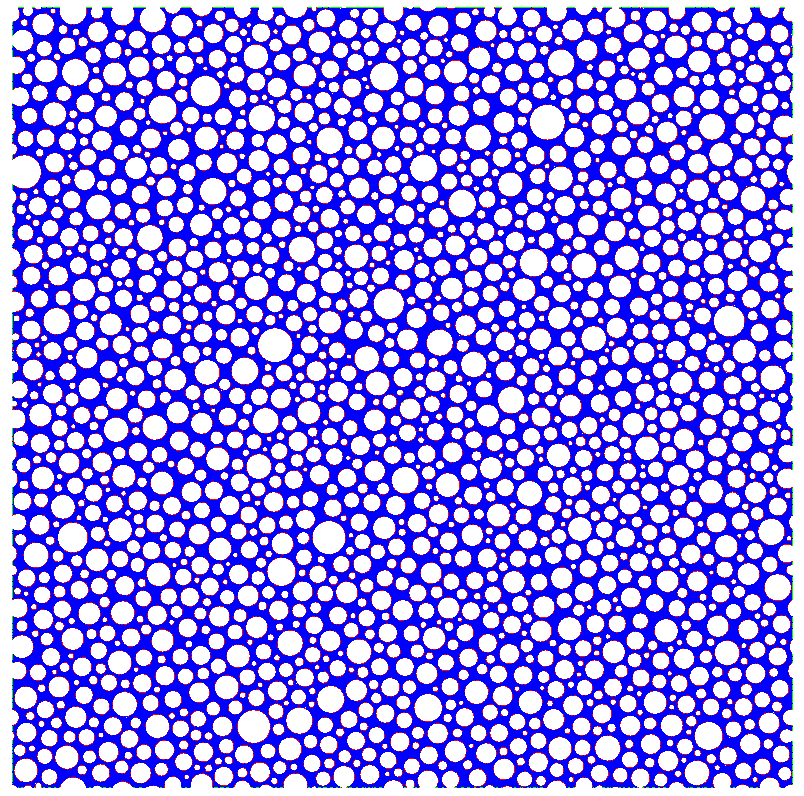}
    \caption{The substrate consists of 1901 non-overlapping circles.}
    \label{fig:substrate1}
\end{figure}

The substrate under consideration consists of 1901 non-overlapping circles which represent the axons and
and we only take into account the regions that are within a square, the domain $\Omega$, that measures $50 \mu m$  
on the side (Figure \ref{fig:substrate1}). The radius of the circles is within the range of 
$0.150 \; \mu m$ up to $1.141 \; \mu m$.
The extracellular region is the exterior of the circles and the diffusion coefficient 
is set to $k_0= 450 \mu m/s^2$.

The \textit{ex vivo} coefficient oscillates between 450 and 600, depending on the temperature and the substances of the medium, \cite{dyrby2013}. We use the smaller one to be able to use small substrates.

It is apparent that numerical approximations are required to simulate the molecular diffusivity is this rather complex
porous medium.

\section{DG-CUDA solution of the heat equation}

\subsection{The Discontinuous Galerkin Method}

\begin{center}
  \begin{figure}[t]    
    \centering  
    \includegraphics[scale=.5]{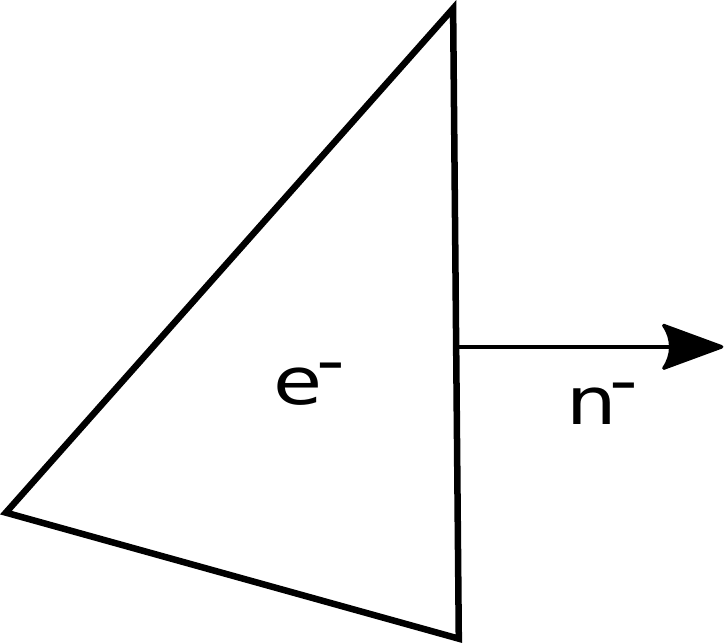}
    \caption{Integration element and its normal.}
      \label{fig:rsol}
  \end{figure}  
\end{center}

Let $\Omega_e$ be partitioned into non overlapping polygonal elements, e.g. a triangulation. Let us denote by $K\equiv K^-$ one of such elements,
see Figure \ref{fig:rsol}.

A defining feature of the DG method is to reduce the PDE to a first order system. Consequently, let us consider

\begin{equation}
  \begin{cases}
    \label{eq:rpde}
    \mathbf{q} &=  \nabla u \\
    \frac{\partial u}{\partial t} &= \nabla \cdot \left(k \mathbf{q}\right).
  \end{cases}
\end{equation}

Multiplying (\ref{eq:rpde}) by the test functions $\mathbf{v}=(v_1,v_2)$, and $v$, we obtain after integrating by parts,

\begin{align*}
  \int_{K^-} \mathbf{q} \cdot \mathbf{v} d\mathbf{x} &= \int_{K^-} ( \nabla u)\cdot\mathbf{v} d\mathbf{x} \\
                  %                           &= \int_{V} \nabla \cdot (  u \mathbf{v}) d\mathbf{x} - \int_{V} u \cdot \nabla \mathbf{v}d\mathbf{x}  \\
                                                   &=  \int_{ \partial K^-} (u \mathbf{v}\cdot \mathbf{n})d\mathbf{s} - \int_{K^-} u \cdot \nabla  \mathbf{v} d\mathbf{x}
\end{align*}

\begin{align*}
  \int_{K^-} u_t v d\mathbf{x} &= \int_{K^-} \nabla  \cdot ( k\mathbf{q})v d\mathbf{x} \\
     %                                        &= \int_{ C} \nabla \cdot ( \sqrt{k} \mathbf{u} \mathbf{v}) d\mathbf{x} - \int_{C} \mathbf{u} \nabla \cdot (\sqrt{k}\mathbf{v})d\mathbf{x}  \\
  &=  \int_{ \partial K^-} v ( k  \mathbf{q})  \cdot \mathbf{n}d\mathbf{s} - \int_{K^-}  k \mathbf{q} \nabla \cdot v d\mathbf{x}
\end{align*}

A second feature of the DG method is the element-wise  approximation of the unknown $\mathbf{q}$ and $u$. Continuity is not enforced at the boundary 
of adjacent elements. Consequently, the boundary terms $u \mathbf{n}$, $k  \mathbf{q}  \cdot \mathbf{n}$ are replaced by boundary
fluxes $h_{u,K^-}(u^-,u^+,\mathbf{n}^-)$, $h_{\mathbf{q},K^-}(\mathbf{q}^-,\mathbf{q}^+,\mathbf{n}^-)$. As customary, the $-$ superscript denotes limits from the interior of $K^-$, and the $+$ superscript  limits from the exterior.

\bigskip

This yields in element $K^-$
\begin{equation}
\begin{array}{rcl}
  \int_{K^-} \mathbf{q} \cdot \mathbf{v} d\mathbf{x} & = & 
  \int_{ \partial K^-}\mathbf{v}\cdot h_{u,K^-}(u^-,u^+,\mathbf{n}^-)d\mathbf{s} - \int_{K^-} u \cdot \nabla  \mathbf{v} d\mathbf{x} \\
   & & \\
   \int_{K^-} u_t v d\mathbf{x} & = &
     \int_{ \partial K^-} v h_{\mathbf{q},K^-}(\mathbf{q}^-,\mathbf{q}^+,\mathbf{n}^-)d\mathbf{s} - \int_{K^-}  k \mathbf{q} \nabla \cdot v d\mathbf{x}
\end{array}
\label{weak_l}
\end{equation}

\bigskip

Let $w$ be an approximation of any of the scalar funcions $q_1$, $q_2$ ,$u$. Within $K^-$ the approximation is given by
\begin{equation}
w(x,t)=\sum_{j=0}^pw_j(t)N_j(x).
\end{equation}

For a triangulation, the functions $N_j$ are chosen as in the Finite Element Method with lagrangian interpolation.

It is assumed that $u^+$ and $\mathbf{q}^+$ are known in (\ref{weak_l}). Hence, we are led to solve an differential-algebraic system for
$u^-$ and $\mathbf{q}^-$. The solution in time is advanced by a Runge-Kutta method.

\bigskip

\noindent\textbf{Remark. }We stress that the solution of the $p+1$ differential-algebraic system (\ref{weak_l}), is solved independently for each element. In practice $p\leq 3$ suffices. Thus, we have small systems to solve that do not exchange information in each time step.

\bigskip

\subsection{Numerical flux}
 There is no  preferred direction of propagation in the heat equation, thus for $u$ a central flux is considered, namely 

\[
h_{u,K^-}(u^-,u^+,\mathbf{n}^-)=\frac{u^-+u^+}{2}\mathbf{n}^-.
\]

A physical assumption is that there is no flow between axons and the extracellular region. Consequently, for $\mathbf{q}$, we propose the numerical flow
\[
h_{\mathbf{q},K^-}(\mathbf{q}^-,\mathbf{q}^+,\mathbf{n}^-)=
\frac{2k^-k^+}{k^-+k^+}\frac{1}{2}(\mathbf{q}^-+\mathbf{q}^+)\cdot\mathbf{n}^-.
\]

This is coined for the problem under consideration. For instance, if $k^->0$ and $k^+=0$, the harmonic mean forces a zero Neumann condition, hence there is no flow trough  the boundary of the element $K^-$ as expected.

\subsection{Cuda implementation}

\begin{center}
  \begin{figure}[t]    
    \centering
    \label{fig:triangulacion}
    \includegraphics[scale=.8]{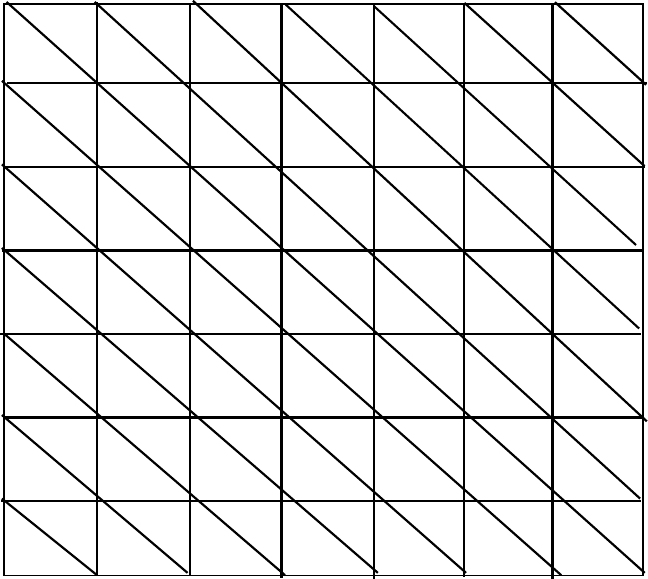}
    \caption{Meshing.}
  \end{figure}  
\end{center}

Meshing is a time consuming task in Galerkin type methods, as it is boundary conditions handling when assembling the local systems. In our case, the domain is divided in square pixels which we use to our advantage. More precisely, these squares are separated by the diagonal in two triangles. A triangle constructed in this fashion, is the basic element for discretization. See Figure \ref{fig:triangulacion}

Also in this MRI application, the heat equation is solved in the extracellular region $\Omega_e$ where  Cauchy data is prescribed. Apparently, there is no need to consider the axon region $\Omega_a$. Nevertheless, we solve the heat equation in elements contained in $\Omega_a$ where
the contribution to the solution is null. 

We are led to balanced computations on every element. As pointed out in the remark above, the calculations in (\ref{weak_l}) are element independent when updating time. Consequently the main ingredients for parallel processing using GPUs are met. Namely, small balanced computations with no exchange of information between processors. 

Parallel processing using GPUs is implemented in the CUDA language in a DELL laptop with hardware:

\noindent CPU: Intel(R) Core(TM) i7-6820HQ CPU @ 2.70 GHz

\noindent GPU: NVIDIA Corporation GM107GLM [Quadro M1000M]

\section{The Gaussian profile of the extra-cellular diffusivity}

In this section we introduce a scheme to approximate the covariance matrix of the EAP. Hence describing the 
Gaussian profile of the extra-cellular diffusivity. Numerical results are also shown for the case of free diffusion, and a comparison with MCMC.

\subsection{The scheme}
Let  us choose $(x_i,y_i)\in\Omega_e$, $i=1,2,\ldots,m$ randomly and uniformly in $\Omega_e$.  Set a final time $T$.
\begin{itemize}
\item For $i=1,2,\ldots,m$, let $U_i(x,y,T)$ be the DG-Cuda solution of the IVBP (\ref{diff_eq})$-$(\ref{Bp_data}), where the Cauchy data is the Dirac's delta function supported in $x_i$.

\item For $i=1,2,\ldots,m$, let $u_i(x,y,T)$ be obtained from $U_i(x,y,T)$ by centering and normalization to yield a density function

\item Construct the mixture model
\[
u(x,y)=\frac{1}{m}\sum_{i=1}^mu_i(x,y,T).
\]

\item Fit a Gaussian density to $u$, that is, determine a covariance matrix $\Sigma$, such that
\begin{equation*}
  u(x,y;\Sigma) \approx \frac{1}{2\pi \sqrt{|\Sigma}|} \exp \left(-\frac{1}{2} \left[ ( (x,y)-(\mu_x,\mu_y))^T \Sigma^{-1} ((x,y)-(\mu_x,\mu_y)) \right] \right),
    \label{eqn:gaussian}
\end{equation*}
\end{itemize}

The covariance matrix is given by,

\begin{displaymath}
  \Sigma =
  \begin{pmatrix}
    E[(EX-\mu_x)(EX-\mu_x)] & E[(EX-\mu_x)(EY-\mu_y)] \\
    E[(EY-\mu_y)(EX-\mu_x)]  & E[(EY-\mu_y)(EY-\mu_y)]
  \end{pmatrix}.
\end{displaymath}

It is computed by quadrature rules using $u_{ij}$, the values of the numerical solution at the nodes $(x_i, y_j)$ in the mesh. 

\subsection{Numerical results}
The substrate under study occupies a square domain of side $50 \;\mu m$.  Cauchy data is given within a centered box of side $20 \;\mu m$.  The observation time is $t=0.036$ seconds. At this given time the outer boundary is not reached by diffusion. In all cases a square grid of $K\times L\equiv 400\times 400$ mesh is used.

\bigskip

\begin{figure}
\centering
  \centering
  \includegraphics[scale=.2]{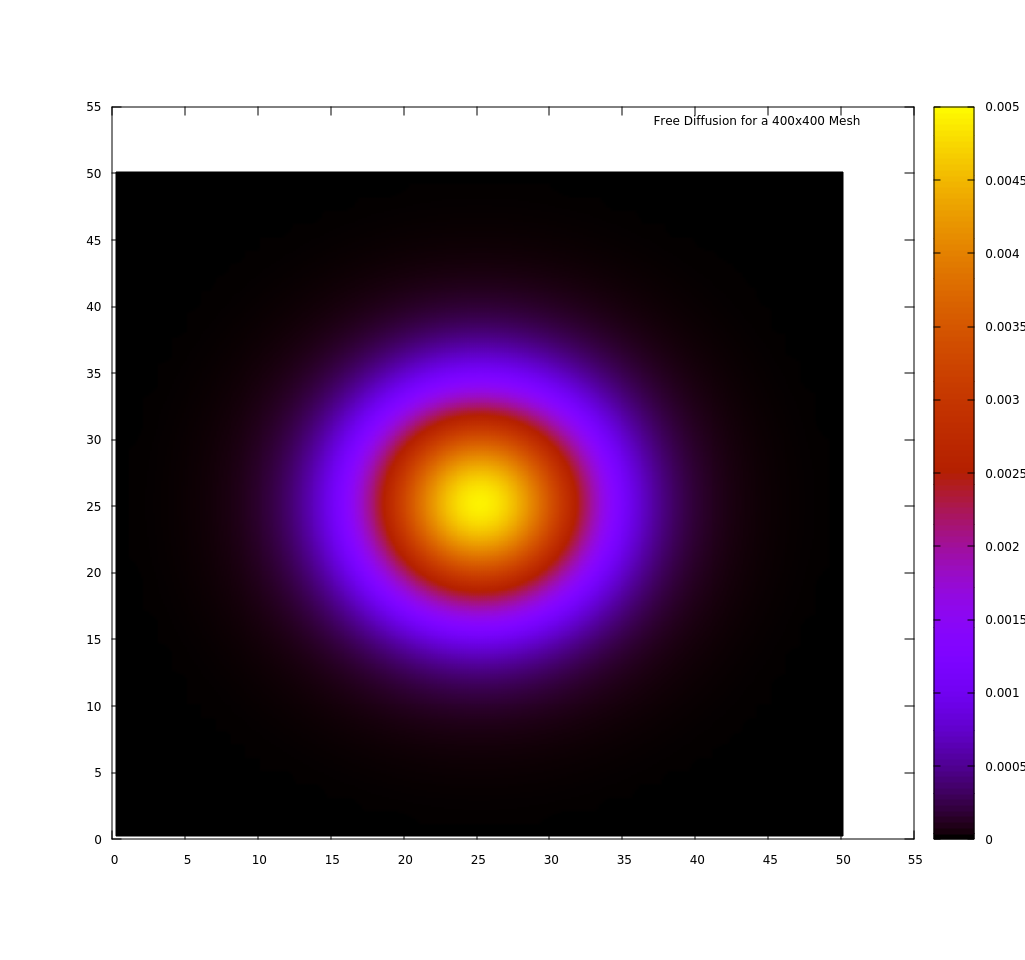}
  \includegraphics[scale=.2]{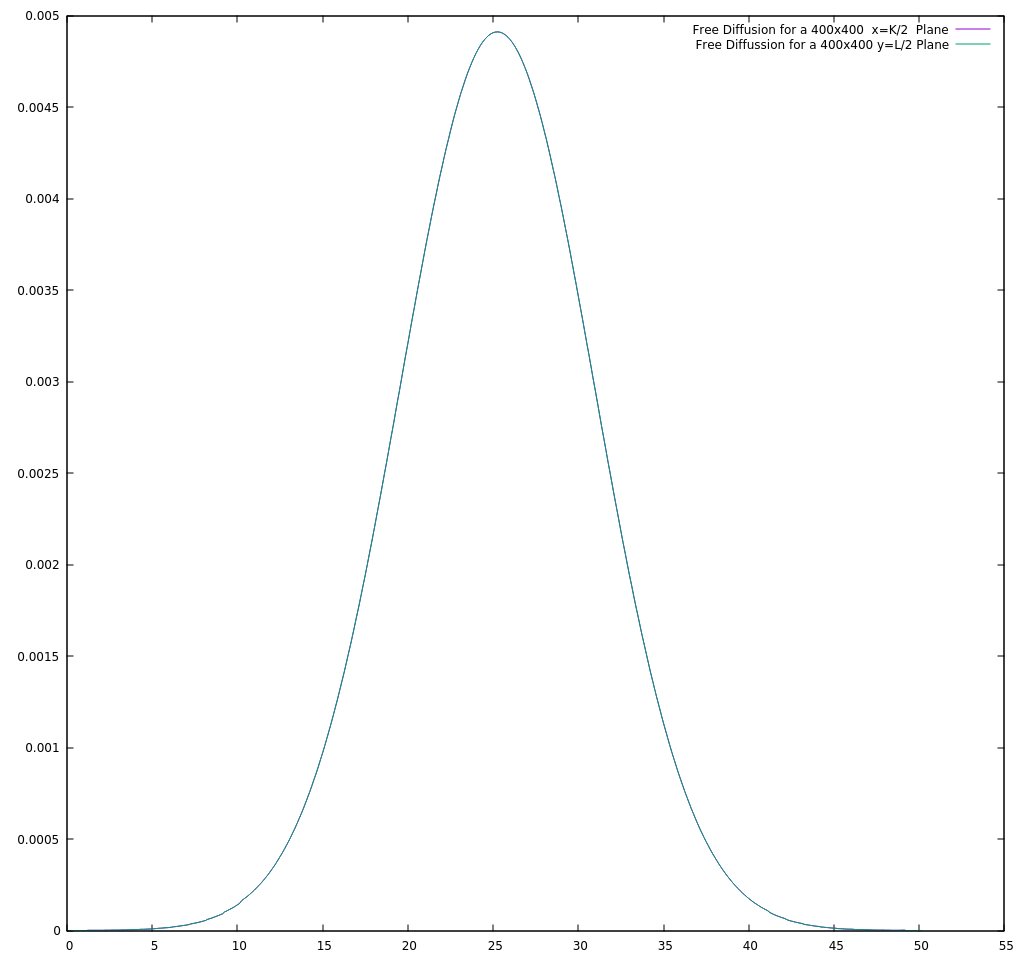}
  \caption{Left: 2D free diffusion. Right: 1D view.  }
  \label{fig:FD2}
\label{fig:FD1}
\end{figure}

\noindent\textbf{Free diffusion. }In this simple example, We solve the Heat Equation with Cauchy data a Dirac's delta supported at the origin. The uniform diffusion coefficient is $k_0= 450 \mu m/s^2$. Graphical results in Figure \ref{fig:FD1}.

Let us compare the free diffusion (without axons) approximation,  versus the
analytical solution. The latter is a bivariate Gaussian function $f(x,y;\Sigma)$ defined in \eqref{eqn:gaussian}, where 
the covariance matrix is 
\begin{displaymath}
  \Sigma =
  \begin{pmatrix}
    2Tk & 0 \\
    0  & 2Tk
  \end{pmatrix}.
\end{displaymath}

Taking $k = 450 \mu m/s^2$ and $T=0.036 s$, non zero coefficients in the covariance
matrix are equal to $32.4$.

\bigskip

The DG-CUDA Gaussian matrix fit for $400\times400$ is
\begin{displaymath}
  \begin{pmatrix}    
    32.53 &    0.00038 \\
   0.00038  &   32.53
  \end{pmatrix}
\end{displaymath}

An alternative construction of the covariance matrix is by means of Monte Carlo Diffusion Simulation. See  \cite{article-MattHall2009}
for details. We just list the corresponding data in their notation.

Diffusion constant $D = 4.5e-10$, $t_s =0.036$ duration of the diffusion simulation. $T=5000$ is the number of time steps in the
simulation. The step length $l$ is obtained from the relation
\[
l=\sqrt{4\frac{Dt_s}{T}}\approx  0.11\mu m.
\]

The obtained MC  Gaussian matrix is
\begin{displaymath}
  \begin{pmatrix}    
     32.487358329 & -0.075889940 \\
     -0.075889940 & 32.378282498
    \end{pmatrix}
\end{displaymath}

To gauge the approximations we compute  the least squares residual

        \begin{equation}
   \sum_i\sum_j [ u(x_i,y_i;\Sigma) - u_{ij} ]^2.
            \label{eqn:lsp}
        \end{equation}

In both cases the approximation of the Gaussian function is highly accurate. The least squares residual (\ref{eqn:lsp}) is of the order
$O(10^{-8})$. For practical purposes, the Gaussian density functions coincide. But the DG-CUDA approximation is structurally more 
consistent. The matrix is symmetric, the values in the diagonal coincide and the other terms are near zero.

\bigskip

\begin{figure}
\centering
  \centering
  \includegraphics[scale=.2]{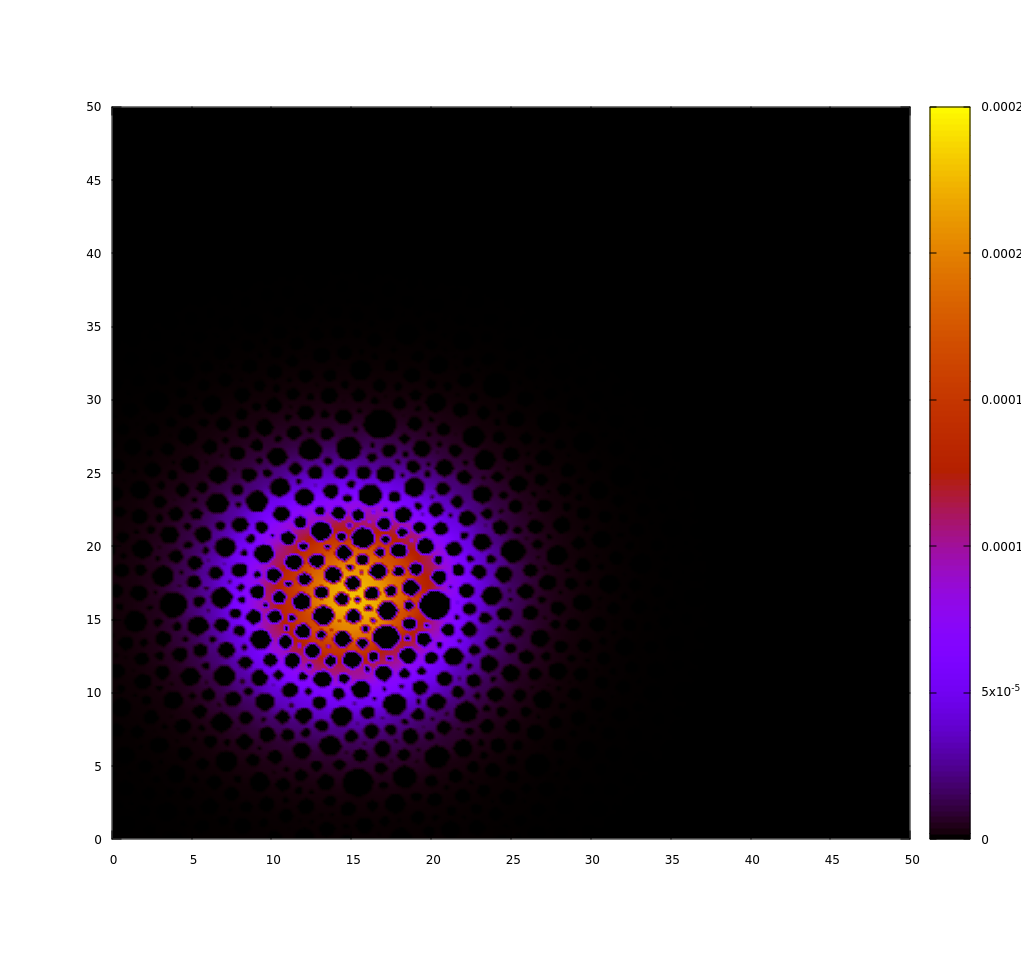}
  \includegraphics[scale=.2]{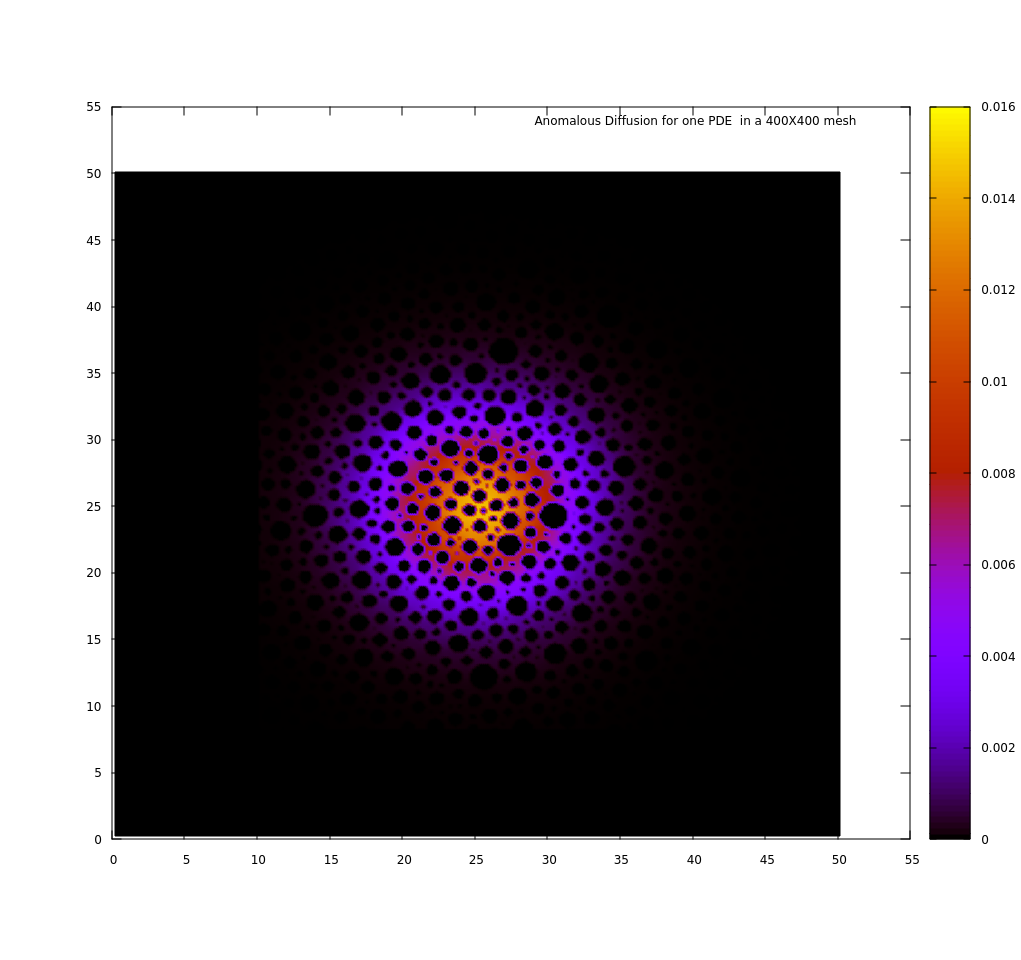}
  \caption{Left: One PDE solution. Right: Normalized and centered solution for one PDE.  }
  \label{fig:L1delta}
\label{fig:R1delta}
\end{figure}

\begin{figure}
\centering
  \centering
  \includegraphics[scale=.2]{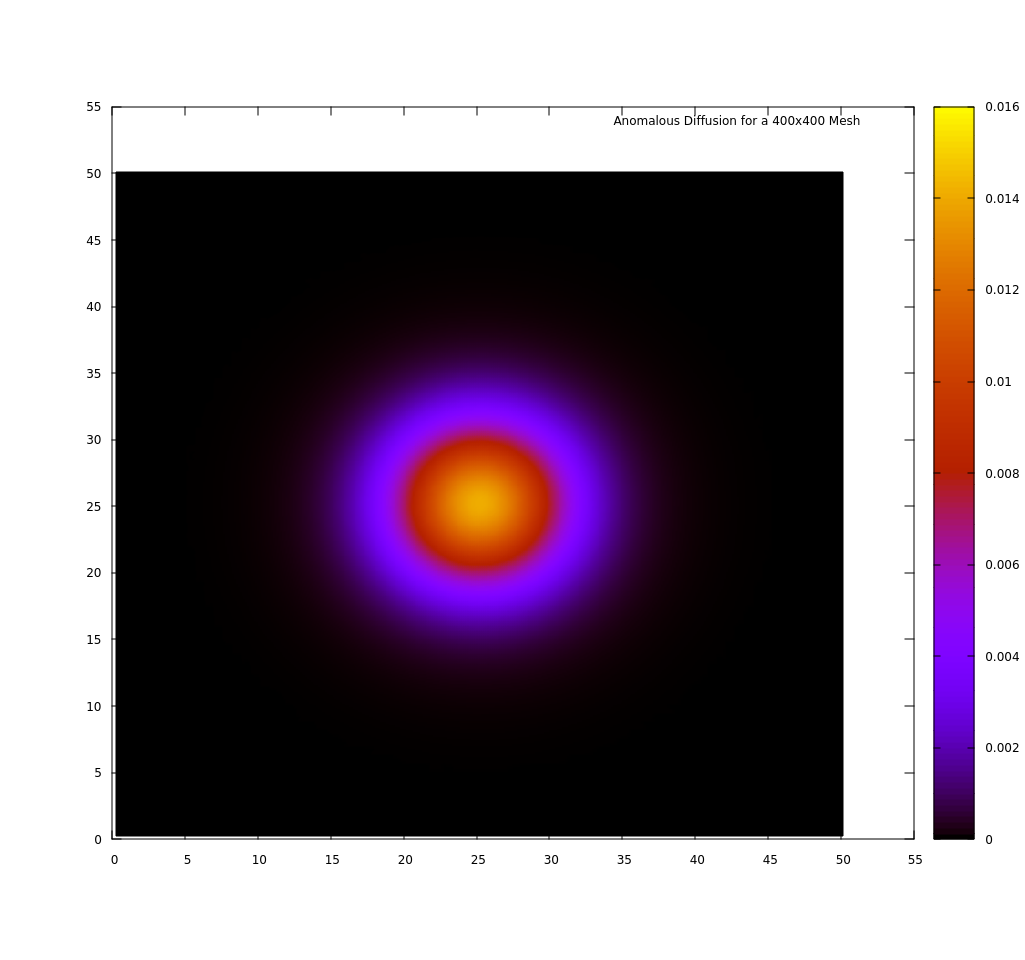}
  \includegraphics[scale=.2]{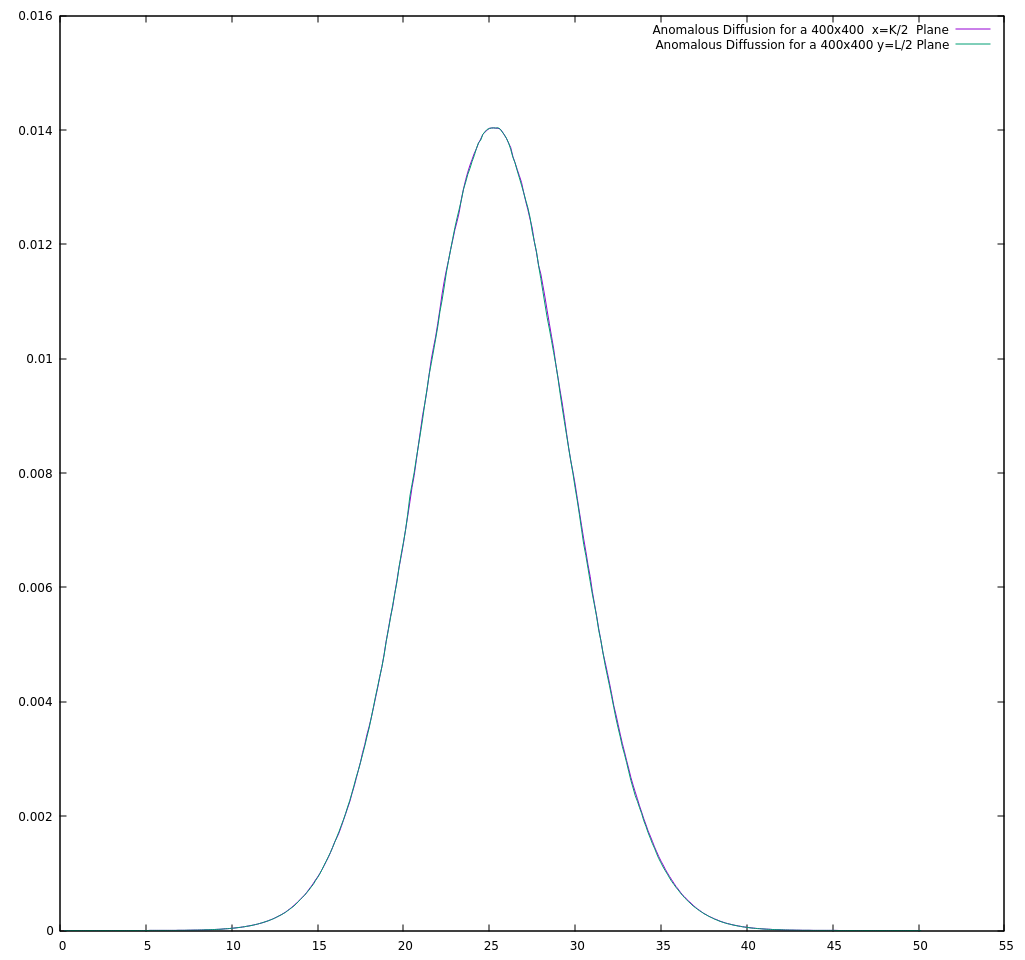}
  \caption{Left: 2D free diffusion. Right: 1D view.  }
  \label{fig:AD2}
\label{fig:AD1}
\end{figure}

\noindent\textbf{Case study. }The axon region $\Omega_a$ is  defined by 1901 axons. The diffusion coefficient as in (\ref{k_coeff}). The scheme above is applied to a mixture of $m=37$ densities. An instance of one PDE solution is shown in Figure \ref{fig:R1delta}. The solution of the full scheme is in Figure \ref{fig:AD1}.
  
\bigskip

The covarince matrix by the DG-CUDA algorithm is

\begin{displaymath}
  \begin{pmatrix}    
    19.50 &    -0.0088 \\
   -0.0088  &   19.50
  \end{pmatrix}
\end{displaymath}

Finally, let us summarize execution time in the following table
\begin{center}
    \begin{tabular}{ | l | c |c| c| c | }
      \hline
      $K\times L$ & Time stps &No. Deltas & CPU  & GPU    \\       \hline \hline
      400$\times$400 & 5184 & 1 & 6585.47 sec. & 5.13649 sec. \\
      400 $\times$400 & 5184 & 37 & 21877.252 sec. & 188.597 sec. \\
      % % \vspace{.2cm}
      \hline
    \end{tabular}
  \end{center}

We remark that regardless of the diffusion coefficient, the GPU process is more than 1000 times faster for the solution of one instance of the heat equation with a single delta as Cauchy data.

\section{Conclusions and future work}

In this work we have provided a methodology for the efficient computation of the diffusivity properties in porous media. These results can be used for the simulation of the DWMR signal given the numerical estimation of the EAP.  Such a computation can be analytically performed by the application of the Fourier transform on the EAP, however, simplistic assumptions about the machine model have to be made \cite{price1997}. In order to provide a useful tool for the medical researchers, the numerical simulator of the signal should take into account the signal changes associated to a realistic MR machine parameters, as the non squared (but trapezoidal) and the finite duration of the magnetic pulses. It is our contention that our methodology is versatile to include this complexities in a 3D extension. The latter is part of our current and future work.

\bigskip

\begin{center}
\textbf{Acknowledgements}
\end{center}

D. Cervantes and M. A. Moreles thank the support of ECOS-NORD through the project: 000000000263116/M15M01. M. A. Moreles also acknowledges  partial support from the project CIB CONACYT 180723.

\bibliographystyle{plain}
\bibliography{Bib_DG_dMRI}

\begin{thebibliography}{10}

\bibitem{article-DistAboitiz1992}
F.~Aboitiz, A.B. Scheibel, R.S.Fisher, and E.Zaidel.
\newblock Fiber composition of the human corpus callosum.
\newblock {\em Brain Res.}, 143(53), 1992.

\bibitem{IntroductionCDMRI}
D.C. Alexander.
\newblock An introduction to computational diffusion mri: the diffusion tensor
  and beyond.
\newblock In Hagen~H. Weickert~J., editor, {\em Visualization and Processing of
  Tensor Fields}, chapter~10, pages 83--106. Springer Verlag, 2006.

\bibitem{article-Alexander2009}
D.C. Alexander.
\newblock Modelling, fitting and sampling in diffusion mri.
\newblock {\em Magnetic Resonance}, 103:255--260, 2009.

\bibitem{article-timeDependenceBurcaw2015}
L.~M. Burcaw, E.~Fieremans, and D.~S. Novikov.
\newblock Mesoscopic structure of neuronal tracts from time-dependent
  diffusion.
\newblock {\em NeuroImage}, 114:18--37, 2015.

\bibitem{cockburn2012discontinuous}
Bernardo Cockburn, George~E Karniadakis, and Chi-Wang Shu.
\newblock {\em Discontinuous Galerkin methods: theory, computation and
  applications}, volume~11.
\newblock Springer Science \& Business Media, 2012.

\bibitem{dyrby2013}
Tim~B Dyrby, Matt~G Hall, Maurice Ptito, Daniel Alexander, et~al.
\newblock Contrast and stability of the axon diameter index from microstructure
  imaging with diffusion {MRI}.
\newblock {\em Magnetic Resonance in Medicine}, 70(3):711--721, 2013.

\bibitem{Nedjati_neuroimage_2017}
Nedjati-Gilani Gemma~L. et. al.
\newblock Machine learning based compartment models with permeability for white
  matter microstructure imaging.
\newblock {\em Neuroimage}, 150:119?135, 2017.

\bibitem{Moroney_MRM24501}
Moroney~B. F., Stait-Gardner T., Ghadirian B., Yadav~N. N., and Price WS.
\newblock Numerical analysis of nmr diffusion measurements in the short
  gradient pulse limit.
\newblock {\em Journal of Magnetic Resonance}, 234:165--175, 2013.

\bibitem{Ferizi2015}
Uran Ferizi, Torben Schneider, Thomas Witzel, Lawrence~L. Wald, Hui Zhang,
  Claudia~A.M. Wheeler-Kingshott, and Daniel~C. Alexander.
\newblock White matter compartment models for in vivo diffusion {MRI} at 300
  mt/m.
\newblock {\em NeuroImage}, 118:468--483, 2015.

\bibitem{Fieremans2010MonteCS}
Els Fieremans, Dmitry~S. Novikov, Jens~H. Jensen, and Joseph~A. Helpern.
\newblock Monte carlo study of a two-compartment exchange model of diffusion.
\newblock {\em NMR in biomedicine}, 23 7:711--24, 2010.

\bibitem{Gilani2017AnIM}
Nima Gilani, Paul~N. Malcolm, and Glyn Johnson.
\newblock An improved model for prostate diffusion incorporating the results of
  monte carlo simulations of diffusion in the cellular compartment.
\newblock {\em NMR in biomedicine}, 30 12, 2017.

\bibitem{article-MattHall2009}
M.G. Hall and D.C. Alexander.
\newblock Convergence and parameter choice for monte-carlo simulations of
  diffusion mri.
\newblock {\em IEEE TRANSACTIONS ON MEDICAL IMAGING}, 28(9), 2009.

\bibitem{DiffusionMRI}
Derek Jones.
\newblock {\em Diffusion MRI theory methods and applications}.
\newblock Oxford University Press, 2011.

\bibitem{article-Neuman1974}
C.H. Neuman.
\newblock Spin echo of spins diffusing in a bounded medium.
\newblock {\em Chemical Physics}, 60(11):4508?4511, 1974.

\bibitem{Nguyen2014AFE}
Dang~Van Nguyen, Jing-Rebecca Li, Denis Grebenkov, and Denis~Le Bihan.
\newblock A finite elements method to solve the bloch-torrey equation applied
  to diffusion magnetic resonance imaging.
\newblock {\em J. Comput. Physics}, 263:283--302, 2014.

\bibitem{article-CompartmentModelsPanagiotaki2012}
E.~Panagiotaki, T.~Schneider, B.~Siow, M.G. Hall, M.F. Lythgoe, and D.C.
  Alexander.
\newblock Compartment models of the diffusion mr signal in brain white matter:
  A taxonomy and comparison.
\newblock {\em NeuroImage}, 2012.

\bibitem{Pierpaoli1996}
C~Pierpaoli, P~Jezzard, P~J Basser, A~Barnett, and G~Di Chiro.
\newblock Diffusion tensor {MR} imaging of the human brain.
\newblock {\em Radiology}, 201(3):637--648, 1996.

\bibitem{price1997}
W.S. Price.
\newblock Pulsed-field gradient nuclear magnetic resonance as a tool for
  studying translational diffusion: Part 1. basic theory.
\newblock {\em \mbox{G-Animal's} Journal}, 1997.

\bibitem{ramirez-manzanaresISMRM2018}
Alonso Ramirez-Manzanares, Mario Ocampo-Pineda, Jonathan Rafael-Patino, Giorgio
  Innocenti, Jean-Philippe Thiran, and Alessandro Daducci.
\newblock Quantifying diameter overestimation of undulating axons from
  synthetic dw-mri.
\newblock In {\em Procc. In Annual Meeting of the International Society of
  Magnetic Resonance in Medicine}, page 748. Annual Meeting of the SMRM,
  París, France., 2018.

\bibitem{article-VanGelderen1994}
P.~van Gelderen, D.~DesPers, C.M. van Zijl, and C.T.W. Moonen.
\newblock Evaluation of restricted diffusion in cylinders. phosphocreatine in
  rabbit leg muscle.
\newblock {\em Magnetic Resonance}, 103:255--260, 1994.

\end{thebibliography}

\end{document}